\documentclass[12pt]{amsart}
\usepackage{amssymb,latexsym,verbatim}
\theoremstyle{plain}
\newtheorem{thm}{Theorem}[section]
\newtheorem{lem}[thm]{Lemma}
\newtheorem{cor}[thm]{Corollary}

\theoremstyle{definition}

\theoremstyle{remark}

\DeclareMathOperator{\Ext}{Ext} \DeclareMathOperator{\Supp}{Supp}
\DeclareMathOperator{\V}{V} 
 
 \DeclareMathOperator{\depth}{depth}
\DeclareMathOperator{\cd}{cd} 
\DeclareMathOperator{\Coass}{Coass} 
\DeclareMathOperator{\Max}{Max} \DeclareMathOperator{\lc}{H}
 
 \DeclareMathOperator{\G}{\Gamma}
 
 \DeclareMathOperator{\pd}{pd}

\newcommand{\lo}{\longrightarrow}
\newcommand{\fa}{\mathfrak{a}}

\newcommand{\fm}{\mathfrak{m}}
\newcommand{\fp}{\mathfrak{p}}

\begin{document}

\title[vanishing and finiteness and Artinianness ]
 {On the vanishing and finiteness properties of
generalized local cohomology modules \\}

\author{Moharram Aghapournahr}
\address{ Moharram Aghapournahr\\   Arak University, Beheshti St, P.O. Box:879, Arak, Iran}
\email{m-aghapour@araku.ac.ir}




\begin{abstract}
 Let $R$ be a commutative noetherian ring, $\fa$ an ideal of $R$ and  $M,N$ finite $R$--modules. We prove
 that the following statements are equivalent.
 \begin{enumerate}
 \item[(i)] $\lc^{i}_{\fa}(M,N)$ is finite for all $i< n$.
 \item[(ii)] $\Coass_R(\lc^{i}_{\fa}(M,N)) \subset \V{(\fa)}$ for all $i< n$.
 \item[(iii)] $\lc^{i}_{\fa}(M,N)$ is coatomic for all $i< n$.
 \end{enumerate}
  If $\pd M$ is finite and $r$ be a non-negative integer such that $r>\pd M$ and $\lc^{i}_{\fa}(M,N)$
is finite (resp. minimax) for all $i\geq r$, then
$\lc^{i}_{\fa}(M,N)$ is zero (resp. artinian)
  for all $i\geq r$.
\end{abstract}

\maketitle

\section{Introduction}
Throughout $R$ is a commutative noetherian ring. Generalized local
cohomology was given in the local case by J. Herzog \cite{He} and in
the more general case by M.H Bijan-Zadeh \cite{BZ}. Let $\fa$ denote
an ideal of a ring $R$. The generalized local cohomology defined
 by
\begin{center}
$\lc^{i}_{\fa}(M,N) \cong \underset{n}\varinjlim
\Ext^{i}_{R}(M/{\fa}^{n}M,N).$
\end{center}

This concept was studied in the articles \cite{S}, \cite{He} and
\cite{Ya}. Note that this is in fact a generalization of the usual
local cohomology, because if $M=R$, then
$\lc^{i}_{\fa}(R,N)=\lc^{i}_{\fa}(N)$. Important problems concerning
local cohomology are vanishing, finiteness and artinianness results
(see \cite{Hu}).

In Section 2 we show in \ref{T:kofinlc} that if $M$ is finite and
all generalized local cohomology modules $\lc^{i}_{\fa}(M,N)$ are
coatomic for all $i<n$, then they are finite for all $i<n$. In fact
this is another condition equivalent to Falting's Local-global
Principle for the finiteness of generalized local cohomology modules
(see \cite[Theorem 2.9]{AK}). In Theorem \ref{T:yo} we generalize
Yoshida's theorem ( \cite[Theorem 3.1]{Yo}).

In Section 3, We prove in \ref{T:mmart}, that when $M$ is a finite
$R$--module of finite projective dimension such that the generalized
local cohomology modules
 $\lc^{i}_{\fa}(M,N)$ are minimax modules for all $i\geq r$, (where $r> \pd M$)
  then they must be artinian.

 For unexplained terminology we refer to \cite{BSh} and \cite{BH}.


\section{Finiteness and vanishing}
An $R$--module $M$ is called {\it coatomic} when each proper
submodule $N$ of $M$ is contained in a maximal submodule
$N^{\prime}$ of $M$ (i.e. such that $M/N^{\prime}\cong R/\fm$ for
some $\fm \in \Max{R}$). This property can also be expressed by
$\Coass_R(M)\subset \Max{R}$ or equivalently that any artinian
homomorphic image of $M$ must have finite length. In particular all
finite modules are coatomic. Coatomic modules have been studied by
Z\"{o}schinger \cite{Zrko}.
\begin{thm}\label{T:kofinlc}
Let $R$ be a noetherian ring, $\fa$ an ideal of $R$ and $M,N$ finite
$R$--modules. The following statements
 are equivalent:
\begin{enumerate}
  \item[(i)] $H^{i}_{\fa}(M,N)$ is coatomic for all $i<n$.
  \item[(ii)] $\Coass_R(\lc^i_\fa(M,N)) \subset \V{(\fa)}$ for all $i<n$.
  \item[(iii)] $H^{i}_{\fa}(M,N)$ is finite for all $i<n$.
\end{enumerate}
\end{thm}
\begin{proof}
By \cite[Theorem 2.9]{AK} and\cite[1.1, Folgerung]{Zrko} we may
assume that $(R,\fm)$ is a local ring.
 \item[(i)]$\Rightarrow$ (ii) It is trivial by the definition of coatomic modules.
 \item[(ii)]$\Rightarrow$ (iii) By \cite[Satz 1.2]{Zrkoass} there is $t \geq 1$ such that $\fa^t\lc^i_\fa(M,N)$ is
  finite for all $i<n$. Therefore there is $s \geq t$ such that $\fa^s\lc^i_\fa(M,N)=0$ for all $i<n$, and apply
  \cite[Theorem 2.9]{AK}.
 \item[(iii)]$\Rightarrow$ (i) Any finite $R$--module is coatomic.
\end{proof}





The following results are generalizations of \cite[Proposition
3.1]{Yo}.

\begin{thm}\label{T:yo}
 Let $(R,\fm)$ be a local ring, $\fa$ be an ideal of  $R$ and $M$ be a finite module of finite projective dimension.
 Let $N$ be a finite module and $r> \pd{M}$. If $\lc^i_\fa(M,N)$ is finite for all $i\geq r$, then
 $\lc^i_\fa(M,N)=0$ for all $i\geq r$.
\end{thm}
\begin{proof}
 We prove by induction on $d=\dim N$. If $d= 0$, By \cite[Theorem
 3.7]{Ya}, it follows that  $\lc^i_\fa(M,N)=0$ for all $i> \pd M +\dim(
 M\otimes_{R}N)$ and so the claim clearly holds for $n=0$.
 Now suppose $d> 0$ and $\lc^i_\fa(M,N)=0$ for all $i> r$. It is
 enough to show  $\lc^r_\fa(M,N)=0$.
 First suppose $\depth_R{N}> 0$. Take $x\in \fm$ which is $N$--regular.
 Then $\dim{ N/{x}N}= d-1$. The exact sequence
\begin{center}
$0\lo N\overset{x}\lo N\lo N/{x}N\lo 0$
\end{center}
induces the exact sequence
\begin{center}
$\lc^{r}_{\fa}(M,N)\overset{x}\lo \lc^{r}_{\fa}(M,N)\lo
\lc^{r}_{\fa}(M,N/{x}N)\lo \lc^{r+1}_{\fa}(M,N)= 0$
\end{center}
It yields that $\lc^{i}_{\fa}(M,N/{x}N)= 0$ for all $i> r$. Hence by
induction hypothesis we get $\lc^{r}_{\fa}(M,N/{x}N)= 0$. Thus we
have $\lc^{r}_{\fa}(M,N)= 0$ by Nakayama's lemma. Next suppose
$\depth_R{N}= 0$. Put $L=\G_{\fm}(N)$. Since $L$ have finite length,
so we have $\dim L=0$ and  therefore $\lc^{i}_{\fa}(M,L)= 0$ for all
$i> \pd M$ . But from the exact sequence
\begin{center}
$0\lo L\lo N\lo N/L\lo 0$
\end{center}
we get the exact sequence
\begin{center}
$...\rightarrow \lc^{i}_{\fa}(M,L)\rightarrow
\lc^{i}_{\fa}(M,N)\rightarrow \lc^{i}_{\fa}(M,N/L)\rightarrow
\lc^{i+1}_{\fa}(M,L)\rightarrow ...$
\end{center}
hence we have $\lc^{i}_{\fa}(M,N)\cong \lc^{i}_{\fa}(M,N/L)$ for all
$i> \pd M$, and we get the required assertion from the first step.
\end{proof}

\begin{thm}\label{T:vanko1}
Let $\fa$ be an ideal of $R$ and $M$ a finite $R$--module of finite
projective dimension. Let $N$ be a finite $R$--module and $r> \pd
M$. The following statements are equivalent:
\begin{enumerate}
 \item[(i)] $\lc^i_\fa(M,N)=0$ for all $i\geq r$.
 \item[(ii)] $\lc^i_\fa(M,N)$ is finite for all $i\geq r$.
 \item[(iii)] $\lc^i_\fa(M,N)$ is coatomic for all $i\geq r$.
\end{enumerate}
\end{thm}
\begin{proof}
$(i)\Rightarrow (ii)\Rightarrow (iii)$ Trivial.
 $(iii)\Rightarrow(i)$ By use of theorem \ref{T:yo} and \cite[1.1,
Folgerung]{Zrko} we may assume that $(R,\fm)$ is a local ring.
 Note that coatomic modules satisfy Nakayama's lemma. So the proof is the same as in theorem \ref{T:yo}.
\end{proof}
In the following corollary $\cd_{\fa}(M,N)$ denote the supremum of
$i$'s such that $\lc^i_\fa(M,N)\neq0$.

\begin{cor}\label{C:vanco1}
Let $\fa$ an ideal of $R$, $M$ a finite $R$--module of finite
projective dimention and $N$ a finite $R$--module. If
$c:=\cd_{\fa}(M,N)> \pd M$, then $\lc^{c}_{\fa}(M,N)$ is not
coatomic in particular is not finite.
\end{cor}

\section{Artinianness}

Recall that a module $M$ is a {\it minimax} module if there is a
finite (i.e. finitely generated) submodule $N$ of $M$ such that the
quotient module $M/N$ is artinian. Thus the class of minimax modules
includes all finite and all artinian modules.
 Moreover, it is closed under taking submodules, quotients and extensions, i.e., it is a Serre subcategory of the
 category of $R$--modules. Minimax modules have been studied by Zink in \cite{Zi} and  Z\"{o}schinger in
  \cite{Zrmm,Zrrad}. See also \cite{Ru}.

\begin{lem}\label{L:flat}
Let $M$ and $N$ be two $R$--module. If $f: R \longrightarrow S$ is a
flat ring homomorphism, then
$$\lc^{i}_{\fa}(M,N)\otimes_R{S}\cong
\lc^{i}_{{\fa}}S(M\otimes_R{S},N\otimes_R{S}).$$

\end{lem}
\begin{proof}
It is easy and we lift it to the reader.

\end{proof}


\begin{thm}\label{T:mmart}
Let $\fa$ an ideal of $R$ and $M$ a finite $R$--module of finite
projective dimension. Let $N$ be a finite $R$--module and $r> \pd
M$. If  $\lc^{i}_{\fa}(M,N)$ is a minimax module for all $i \geq r$,
then $\lc^{i}_{\fa}(M,N)$ is an artinian module for all
 $i \geq r$.
\end{thm}
\begin{proof}

Let $\fp$ be a non-maximal prime ideal of $R$. Then by the
definition of minimax module and  lemma \ref{L:flat}
$\lc^{i}_{\fa}(M,N)_\fp \cong \lc^{i}_{{\fa}R_{\fp}}(M_\fp,N_\fp)$
is a finite $R_\fp$--module for all $i\geq r$. By theorem
\ref{T:yo}, $\lc^{i}_{\fa}(M,N)_\fp=0$ for all $i\geq r$,
 thus $\Supp_R(\lc^{i}_{\fa}(M,N))\subset{\Max{R}}$ for all $i\geq r$. By \cite[Theorem 2.1]{Ru},   $\lc^{i}_{\fa}(M,N)$
 is artinian for all $i \geq r$.
\end{proof}
Let $q_{\fa}(M,N)$ denote the supremum of the $i$'s such that
$\lc^{i}_{\fa}(M,N)$ is not artinian with the usual convention that
the supremum of the empty set of integers is interpreted as
$-\infty$.

\begin{cor}\label{C:mmart}
Let $\fa$ an ideal of $R$, $M$ a finite $R$--module of finite
projective dimension and $N$ a finite $R$--module. If
$q:=q_{\fa}(M,N)> \pd M$, then $\lc^{q}_{\fa}(M,N)$ is not minimax
in particular is not finite.
\end{cor}

\providecommand{\bysame}{\leavevmode\hbox
to3em{\hrulefill}\thinspace}

\end{document}